% Paper "Continuous families of hyperfinite subfactors with
% the same standard invariant"
%
% Dietmar Bisch, Remus Nicoara, Sorin Popa
%
% This is a plain TeX file, use "tex " to compile
% (amstex 2.1 and epsf.sty are used via the \input command)
%
%
\input amstex
\documentstyle{amsppt}
\def\twoone{${\hbox{\uppercase\expandafter
{\romannumeral2}}}_1$}

\def\Ad{{\hbox{Ad\,}}}
\def\Aut{{\hbox{Aut}}}
\def\Out{{\hbox{\text{\rm Out}}}}

\def\Int{{\hbox{\rm Int}}}

\TagsOnRight

\baselineskip=18 pt plus 2pt
\lineskip=9 pt minus 1 pt
\lineskiplimit=9 pt

\magnification = \magstep 1
\vsize=7.5in

\topmatter
\title
CONTINUOUS FAMILIES OF HYPERFINITE SUBFACTORS 
WITH THE SAME STANDARD INVARIANT
\endtitle
\rightheadtext{A continuous family of non-isomorphic subfactors}
\author DIETMAR BISCH$^1$, REMUS NICOARA$^2$ AND SORIN POPA$^3$
\endauthor
\leftheadtext{Dietmar Bisch, Remus Nicoara, Sorin Popa}
\thanks 
$^1$ supported by NSF under Grant No. DMS-0301173,
$^2$ supported by NSF under Grant No. DMS 0500933,
$^3$ supported by NSF under Grant No. DMS-0100883
\endthanks
\affil $^{1,2}$ Vanderbilt University \\
Department of Mathematics \\
Nashville, TN 37240 \\
  \\
$^3$ UCLA Mathematics Department\\
Box 951555 \\
Los Angeles, CA 90095-1555
\endaffil

\address 
\email dietmar.bisch\@vanderbilt.edu, remus.nicoara\@vanderbilt.edu, 
\newline
 popa\@math.ucla.edu
\endemail 
\endaddress

\address Vanderbilt University, Department of Mathematics, 
Nashville, TN 37240, USA 
\newline
UCLA Mathematics Department, Box 951555, Los Angeles, CA 90095-1555, USA.
\endaddress


\abstract
We construct numerous continuous families of irreducible subfactors
of the hyperfinite II$_1$ factor, which are non-isomorphic, but have
all the same standard invariant. In particular, we obtain
1-parameter families of irreducible, non-isomorphic subfactors of the
hyperfinite II$_1$ factor, which have Jones index $6$
and have all the same standard invariant with property (T). We exploit 
the fact that property (T) groups have uncountably many non-cocycle 
conjugate cocycle actions on the hyperfinite II$_1$ factor.
\endabstract

\endtopmatter

\document

\heading
0. Introduction
\endheading
The standard invariant $\Cal G_{N,M}$ of an inclusion of 
II$_1$ factors $N \subset M$ with finite
Jones index is an extremely powerful invariant which leads to
a complete classification of all subfactors of the hyperfinite 
II$_1$ factor $R$ with index $\le 4$ (see for instance \cite{Jo4}, 
\cite{GHJ}, \cite{EK}, \cite{Po4}). It turns out that there are 
countably many non-isomorphic subfactors with index $\le 4$, and their
(countably many distinct) standard invariants are enough to 
reconstruct these subfactors. However, when
the Jones index becomes $>4$ the standard invariant will no longer be a
complete invariant for the subfactor in general. In \cite{Po4} a 
notion of {\it amenability} for $\Cal G_{N,M}$ was introduced, 
and it was shown that subfactors $N$ of the hyperfinite II$_1$ factor $R$ 
with amenable $\Cal G_{N,R}$ are classified
by this invariant. It is still open whether a converse of this result
is true. In other words, it is not known whether, given a subfactor 
$P \subset R$ with non-amenable standard invariant $\Cal G_{P,R}$, 
there is another subfactor $Q \subset R$ such that $\Cal G_{P,R}$ and
$\Cal G_{Q,R}$ coincide, but the inclusions $P \subset R$ and 
$Q \subset R$ are not isomorphic. Compare this with the
results of Ocneanu and Jones on outer actions of groups on the
hyperfinite II$_1$ factor (\cite{Oc}, \cite{Jo2}). Namely, it is
shown in \cite{Oc} that an amenable group has only one outer
action on the hyperfinite II$_1$ factor (up to outer conjugacy)
whereas the result in \cite{Jo2} shows that non-amenable groups
have always at least two.

We show in this paper that one can construct uncountably many
examples of irreducible, hyperfinite subfactors with integer 
index which are not isomorphic, but have all the same standard 
invariant. The smallest Jones index for which our construction
works is $6$. This is a rather surprising result since the
standard invariant has so far been sufficiently powerful to 
classify subfactors with small index. Our work shows that $6$ 
has to be considered as a "big" index from this point of view. 
The construction of our exotic subfactors relies mainly on two
ingredients. One is the class of subfactors introduced in \cite{BH}.
Those subfactors are simple quantum dynamical systems that arise
from outer actions of finite groups $H$ and $K$ on a II$_1$ factor
$M$. The subfactor $M^H$ is the fixed point algebra under the $H$
action. It is contained in the crossed product algebra 
$M \rtimes K$. The second ingredient is a rigidity result in \cite{Po6},
which says that infinite discrete groups with Kazhdan's property (T) 
have continuously many non-cocycle conjugate cocycle actions on the 
hyperfinite II$_1$ factor. Since there are many property (T) groups
which can be written as a quotient of $\Bbb Z_2 * \Bbb Z_3 =$ 
PSL($2,\Bbb Z$) (for instance
SL($2n+1,\Bbb Z$) for $n \ge 14$), we obtain irreducible, hyperfinite
subfactors with index $6$ of the form 
$R^{\Bbb Z_2} \subset R \rtimes \Bbb Z_3$,
whose relative fundamental group is trivial. This means that
the subfactors $pR^{\Bbb Z_2}p \subset p(R \rtimes \Bbb Z_3)p$,
where $p \in R^{\Bbb Z_2}$ is a projection of trace $t$, are mutually 
non-isomorphic as $t$ runs through $(0,1]$. They have of course
all the same standard invariant.

Here is a more detailed description of the sections in this article.
In section 1 we collect and prove several results on cocycle actions of
discrete groups on II$_1$ factors. In particular we identify 
the reduction of a crossed product subfactor by a projection with 
a cocycle crossed product built on the reduced factor. We give several
explicit examples of groups with non-cocycle conjugate (cocycle) actions
on the hyperfinite II$_1$ factor. In section 2 we prove the main
result (theorem 2.2). We show that if $G$ is a discrete ICC group,
generated by a finite abelian group $H$ and a cyclic group $K$ 
of prime order, with an outer and ergodic action $\alpha$ on the
hyperfinite II$_1$ factor $R$, then the relative fundamental
group of the subfactor $R^H \subset R \rtimes_\alpha K$ is
contained in the fundamental group of the action $\alpha$. As a
corollary we obtain numerous 1-parameter families of irreducible,
hyperfinite subfactors of index $6$ which have all the same
standard invariant. 

{\it Acknowledgements.} Part of this work was done while D.B. and
R.N. were visiting UCLA. We would like to thank UCLA for its
hospitality.

\vfill
\eject

\heading
1. Preliminaries
\endheading
\bigskip
For the convenience of the reader we collect in this section several 
results on crossed products by cocycle actions. Most of these are 
well-known to experts. 

\definition {1.1 Cocycle Actions} Let $G$ be a discrete group and 
$M$ a II$_1$ factor. Let Aut($M$), $\Cal U(M)$ denote the automorphism 
group, respectively the unitary group of $M$. A {\it cocycle action} 
$\alpha $ of $G$ on $M$ 
is a map $\alpha:G\rightarrow \text{Aut}(M)$ such that there exists 
a map $v:G\times G\rightarrow \Cal{U}(M)$, with the properties:

\ (i) $\alpha_e=$ id and $\alpha_g\alpha_h=\Ad v_{g,h}\alpha_{gh}$, 
for all $g$, $h\in G$,

(ii) $v_{g,h}v_{gh,k}=\alpha_g(v_{h,k})v_{g,hk}$, for all $g$, $h$,
$k\in G$.

\enddefinition

The map $v$ is called a {\it 2-cocycle} for $\alpha$. $v$ is 
{\it normalized} if $v_{g,e}=v_{e,g}=1$ for all $g\in G$, where
$e$ denotes the identity of $G$. 
Any 2-cocycle $v$ can be normalized by replacing it, if necessary, 
by $v'_{g,h}=v^*_{e,e}v_{g,h}$, $g$, $h\in G$ (note that $v_{e,e}$ 
is a scalar since $M$ is a factor).

All 2-cocycles considered from now on will be normalized. All (cocycle) 
actions considered in this paper will be assumed {\it properly outer}, 
i.e. $\alpha_g$ cannot be implemented by unitary elements in $M$, 
for all $g\not= e$. Also, we will usually denote a cocycle action
as a pair $(\alpha, v)$. 

The next lemma shows that the cocycle $v$ is unique up to a perturbation 
by a {\it scalar 2-cocycle} $\mu$. 

\proclaim{Lemma 1.2} If $v$, $v'$ are normalized 2-cocycles for the 
cocycle action $\alpha$ of $G$ on $M$ then $v=\mu v'$ for some 
normalized scalar 2-cocycle $\mu$ (i.e. 
$\mu:G\times G\rightarrow \Bbb T$ satisfying
$\mu_{e,e}=1$ and $\mu_{g,h}\mu_{gh,k}=\mu_{h,k}\mu_{g,hk}$, for all 
$g$, $h$, $k\in G$). 
\endproclaim

\demo{Proof }
 $\Ad v_{g,h}=\Ad v'_{g,h}$, for all $g$, $h\in G$ 
implies $v^*_{g,h}v'_{g,h}\in \Cal Z(M)=\Bbb C$, so there 
exists $\mu:G\times G\rightarrow \Bbb T$ such that $v=\mu v'$. 
Since $v$, $v'$ are normalized we have $\mu_{e,e}=1$. Using 
$v'=\mu v$ in the relation $v'_{g,h}v'_{gh,k}=\alpha_g(v'_{h,k})v'_{g,hk}$ 
it follows 
$$\mu_{g,h}\mu_{gh,k}v_{g,h}v_{gh,k}=\mu_{h,k}\mu_{g,hk}
\alpha_g(v_{h,k})v_{g,hk}$$ 
so $\mu$ satisfies the 2-cocycle relation 
$\mu_{g,h}\mu_{gh,k}=\mu_{h,k}\mu_{g,hk}$, for all $g$, $h$, $k\in G$.
\qed
\enddemo

A 2-cocycle $v$ for the action $\alpha$ is called a {\it coboundary} 
(or a {\it trivial cocycle}) if there exists a map 
$w:G\rightarrow \Cal U(M)$ such that $w_e=1$ and 
$v_{g,h}=\alpha_g(w_h^*)w_g^*w_{gh}$, for all $g$, $h\in G$. 

\definition{1.3 Conjugacy of actions} We say that two cocycle actions 
$(\alpha^1,v^1)$, $(\alpha^2,v^2)$ of the groups $G_1$ resp. $G_2$ on 
the II$_1$ factors $M_1$, 
$M_2$ are {\it cocycle conjugate} if there exists a *-isomorphism 
$\Phi: M_1 \to M_2$ (onto), a group isomorphism 
$\gamma:G_1\rightarrow G_2$ and $w_g\in \Cal U(M_2)$ such that:

\ (i) $\Phi \alpha^1_g \Phi^{-1}=\Ad w_g \circ \alpha^2_{\gamma(g)}$, 
for all $g\in G_1$,

(ii) $\Phi(v^1_{g,h})=w_g \alpha^2_{\gamma(g)}(w_h)v^2_{\gamma(g),
\gamma(h)}w^*_{gh}$, for all $g$, $h\in G_1$.
\enddefinition

The cocycle actions $(\alpha^1,v^1)$, $(\alpha^2,v^2)$ of
$G_1$ resp. $G_2$ are said to 
be {\it outer conjugate} (or {\it weakly cocycle conjugate}) if 
condition (i) holds. If $\alpha^1$, $\alpha^2$ are properly outer, 
(i) is equivalent to saying that 
the images of $G_i$ under $\alpha^i$ in $\Out(M_i) \overset{\text{def}}\to
= \Aut (M_i)/\Int (M_i)$, $i=1$, $2$, are conjugate by a
*-isomorphism $\Phi:M_1\to M_2$. 

Indeed, if $\Phi \alpha^1(G_1) \Phi^{-1}=\alpha^2(G_2)$ 
in Out($M_2$), there exists a bijection $\gamma:G_1\rightarrow G_2$ 
and unitaries $w_g\in \Cal U(M_2)$, such that 
$\Phi \alpha^1_g \Phi^{-1}=\Ad w_g \circ \alpha^2_{\gamma(g)}$, 
for all $g\in G_1$. Since 
$g\rightarrow \Phi \alpha^1_g \Phi^{-1}=\alpha^2_{\gamma(g)} \in$ 
Out($M_2$) is a group morphism and $\alpha^2$ is properly outer, 
it follows that $\gamma$ is a group morphism.

The cocycle actions $\alpha^1$, $\alpha^2$ are called {\it conjugate} 
if both conditions (i), (ii) are satisfied with $w=1$. 

Jones proved that any two outer actions of a finite group on the
hyperfinite II$_1$ factor $R$ are conjugate (\cite{Jo1}). In fact,
any two outer actions of an amenable group on $R$ are cocycle
conjugate (\cite{Oc}). The situation is very different when the
group is not amenable. Any non-amenable group has at least two outer
actions on $R$ which are not outer conjugate (\cite{Jo2}). If the group
is rigid, a much stronger result is true. Following \cite{Po6} we
call a group $G$ {\it weakly rigid} (or {\it w-rigid}) if it has an 
infinite normal
subgroup such that the pair $(G,H)$ has the Kazhdan-Margulis relative 
property (T) (\cite{K}, \cite{Ma}). It is shown in \cite{Po6} 
that if $G$ is w-rigid, then there exists a 
continuous family of non-outer conjugate cocycle actions of $G$ on 
$R$. We will use this fact in the next section.

The next lemma shows that perturbing a cocycle action $\alpha$ by 
unitaries of $M$ gives a cocycle action that is cocycle conjugate 
to $\alpha$.

\proclaim{Lemma 1.4} Let $(\alpha,v)$ be a cocycle action of $G$ on $M$ 
and let $w_g$ be unitaries in $M$ for all $g\in G$. Then 
$\beta_g=\Ad w_g \alpha_g$ 
is a cocycle action of $G$ on $M$ with cocycle $v'$, where 
$v'_{g,h}=w_g\alpha_g(w_h)v_{g,h}w^*_{gh}$ for all $g$, $h\in G$. 
Note that $\beta$ is (trivially) cocycle conjugate to $\alpha$.
\endproclaim

\demo{Proof} We show that $(\beta,v')$ satisfies conditions 1.1 (i) and 
1.1 (ii). We compute 
\medskip
$\beta_g\beta_h=\Ad w_g\alpha_g\Ad w_h\alpha_h$

\noindent\hskip 0.8in 
$=\Ad(w_g\alpha_g(w_h))\alpha_g\alpha_h$

\noindent\hskip 0.8in 
$=\Ad (w_g\alpha_g(w_h))\Ad v_{g,h}\alpha_{gh}$ 

\noindent\hskip 0.8in 
$= \Ad (w_g\alpha_g(w_h)v_{g,h}w_{gh}^*)\Ad w_{gh}\alpha_{gh} $

\noindent \hskip 0.8in 
$ =\Ad (v'_{g,h})\beta_{gh}$

\medskip
\noindent 
for all $g$, $h\in G$, which proves 1.1 (i). We check 1.1 (ii). 

\medskip
\noindent
$\beta_g(v'_{h,k})v'_{g,hk}=(\Ad w_g\alpha_g)(v'_{h,k})v'_{g,hk}$

\noindent\hskip .8in 
$=w_g\alpha_g(w_h\alpha_h(w_k)v_{h,k}w^*_{hk})w_g^*w_g
\alpha_g(w_{hk})v_{g,hk}w^*_{ghk}$

\noindent\hskip .8in 
$=w_g\alpha_g(w_h)\Ad v_{g,h}(\alpha_{gh}(w_k))\alpha_g(v_{h,k})
\alpha_g(w^*_{hk})\alpha_g(w_{hk})v_{g,hk}w^*_{ghk}$

\noindent\hskip .8in 
$=w_g\alpha_g(w_h)v_{g,h}\alpha_{gh}(w_k)(v_{g,h}^*
\alpha_g(v_{h,k})v_{g,hk})w^*_{ghk}$

\noindent\hskip .8in 
$=w_g\alpha_g(w_h)v_{g,h}w^*_{gh}w_{gh}\alpha_{gh}(w_k)v_{gh,k}w^*_{ghk}$

\noindent\hskip .8in  
$=v'_{g,h}v'_{gh,k}$.
\ \  \qed
\enddemo

\medskip
The next lemma shows that if two cocycle actions 
$(\alpha^1,v^1)$, $(\alpha^2,v^2)$ are outer conjugate then there exists 
a scalar 2-cocycle $\mu$ such that $(\alpha^1,v^1)$, $(\alpha^2,\mu v^2)$ 
are cocycle conjugate.

\medskip
\proclaim{Lemma 1.5} If the cocycle actions $(\alpha^1,v^1)$,
$(\alpha^2,v^2)$ of $G_1$, $G_2$ on $M_1$, $M_2$ are outer conjugate by 
$\Phi:M_1 \to M_2$, then there exists a group isomorphism 
$\gamma:G_1\rightarrow G_2$, unitaries $w_g\in \Cal U(M_2)$, for all 
$g\in G_1$ and a scalar 2-cocycle $\mu:G_2\times G_2\rightarrow \Bbb T$ 
such that

\ (i) $\Phi \alpha^1_g \Phi^{-1}=\Ad w_g \circ \alpha^2_{\gamma(g)}$, 
for all $g\in G_1$.

(ii) $\Phi(v^1_{g,h})=\mu_{\gamma(g),\gamma(h)} w_g 
\alpha^2_{\gamma(g)}(w_h)v^2_{\gamma(g),\gamma(h)}w^*_{gh}$, 
for all $g$, $h\in G_1$.
\endproclaim

\demo {Proof} Since $\alpha^1$, $\alpha^2$ are outer conjugate, there 
exist $w_g\in \Cal U(M_2)$, $g\in G_1$, and an isomorphism 
$\gamma:G_1\rightarrow G_2$ such that 
$\Phi\alpha^1_g\Phi^{-1}=\Ad w_g\text{ }\alpha^2_{\gamma(g)}$, 
for all $g\in G_1$.

Since $v^1$ is a 2-cocycle for $\alpha^1$, $\Phi(v^1)$ 
is a 2-cocycle for the cocycle action $\Phi\alpha^1\Phi^{-1}$.  On the 
other hand, Lemma 1.4 implies that 
$v'_{g,h}=w_g \alpha^2_{\gamma(g)}(w_h)v^2_{\gamma(g),\gamma(h)}w^*_{gh}$ 
is a 2-cocycle for the cocycle action 
$g\rightarrow\Ad w_g \circ \alpha^2_{\gamma(g)}=\Phi \alpha^1_g \Phi^{-1}$. 

Since by lemma 1.2 any two 2-cocycles of the 
same cocycle action differ by a scalar 2-cocycle, there exists 
$\mu$ such that 
$\Phi(v^1_{g,h})=\mu_{\gamma(g),\gamma(h)} w_g 
\alpha^2_{\gamma(g)}(w_h)v^2_{\gamma(g),\gamma(h)}w^*_{gh}$, 
for all $g$, $h\in G_1$.  \qed
\enddemo

\medskip
\definition{1.6. Crossed Products by Cocycle Actions} Let $M$ be a 
II$_1$ factor and $\tau$ its unique normalized faithful trace. Let 
$(\alpha,v)$ be a cocycle action of the discrete group $G$ on $M$. 

The {\it crossed product algebra} $(M\rtimes_{\alpha, v} G,\tau)$ is 
defined as the von Neumann subalgebra of $B(l^2(G,L^2(M,\tau)))$ 
generated by unitaries 
$u_g$ $\in B(l^2(G,L^2(M,\tau)))$, $g\in G$, where 
$u_g(f)(h)=v_{g,g^{-1}h}f(g^{-1}h)$, for all $f\in l^2(G,L^2(M,\tau))$, 
$g$, $h\in G$, and by a copy of the 
algebra $M$ given by 
$(x\cdot f)(g)=\alpha_g^{-1}(x)f(g)$, for all $x\in M$, 
$f\in l^2(G,L^2(M,\tau))$, $g\in G$. In this paper we will
most of the time drop the cocycle from the notation and simply write
$(M\rtimes_{\alpha} G,\tau)$.
The formula
$\tau(X)=\langle X\delta_e,\delta_e \rangle$, for all 
$X\in M\rtimes_{\alpha} G$, where $\delta_e\in l^2(G,L^2(M,\tau))$ is 
the $L^2(M,\tau)$-valued function on $G$ that takes value $1$ at $e$ 
and $0$ elsewhere, defines a trace $\tau$ on $M\rtimes_{\alpha} G$. 
See for instance \cite{Su1}, \cite{Su2} for more details.

Alternatively, 
$(M\rtimes_{\alpha} G,\tau)$ can be viewed in the following way: Consider 
the Hilbert algebra $\Cal{M}$ of finite formal sums 
$\Cal M=\{\sum_{g\in G} x_gu_g, x_g\in M\}$, with multiplication rules 
$$u_gu_h=v_{g,h}u_{gh}, \quad u_gx=\alpha_g(x)u_g, \quad x=xu_e=1x,$$ 
for all $g$, $h\in G$, $x\in M$, and $*$-operation
$(u_g x)^*=u_{g^{-1}}\alpha_g(x^*)$. The trace is given by 
$\tau(\sum_{g\in G} x_gu_g)=\tau(x_e)$. 
Then $M\rtimes_\alpha G$ is defined as the closure of $\Cal{M}$ in 
norm $\Vert\text{ } \Vert_{2,\tau}$ on bounded sequences. 
\enddefinition

$(M\rtimes_{\alpha} G,\tau)$ is a finite von Neumann algebra with 
normal faitful trace $\tau$. If the cocycle action $\alpha$ is outer then 
$M'\cap M\rtimes_{\alpha} G=\Bbb C$. In particular, if $\alpha$ is 
outer then $(M\rtimes_{\alpha} G,\tau)$ is a II$_1$ factor.

For the convenience of the reader we include a proof of the well-known
result that the isomorphism class of the {\it inclusion} 
$(M\subset M\rtimes_{\alpha,v} G)$ is determined by the cocycle 
conjugacy class of the (cocycle) action $(\alpha,v)$ (\cite{Jo1}).

\proclaim{Proposition 1.7} Let $\alpha^1,\alpha^2$ be cocycle actions of 
the discrete groups $G_1$, $G_2$ on the II$_1$ factors $M_1$, $M_2$, 
with 2-cocycles $v^1$, $v^2$. If there exists a surjective $*$-isomorphism 
$\Phi: M_1\rtimes_{\alpha^1} G_1 \to  M_2\rtimes_{\alpha^2} G_2$
such that $\Phi(M_1)=M_2$, then
$\alpha^1$ and $\alpha^2$ are cocycle conjugate. More precisely, 
there exists a group isomorphism $\gamma:G_1\rightarrow G_2$, and 
unitaries $w_g\in \Cal{U}(M_2)$, for all $g\in G_1$, such that:

\ (i) $\Phi \alpha^1_g \Phi^{-1}=\Ad w_g\text{ }\alpha^2_{\gamma(g)}$, for 
all $g\in G_1$, 

(ii) $\Phi(v^1_{g,h})=w_g \alpha^2_{\gamma(g)}(w_h)v^2_{\gamma(g),
\gamma(h)}w^*_{gh}$, for all $g,h\in G_1$.

Conversely, if $\Phi:M_1\to M_2$ is a $*$-isomorphism (onto), 
$\gamma:G_1\rightarrow G_2$ is a group isomorphism, and there exist 
unitaries $w_g\in \Cal{U}(M_2)$ for all $g\in G_1$ such that (i), (ii) 
are satisfied, then $\Phi$ can be extended to an isomorphism 
$M_1\rtimes_{\alpha^1} G_1\simeq M_2\rtimes_{\alpha^2} G_2$ (hence
$\Phi$ is an isomorphism of the associated inclusions).

\endproclaim

\demo{Proof} For $i=1$, $2$ let $u_g^i$ denote the unitaries implementing 
the action $\alpha^i$ on $M_i$, i.e. 
$\alpha^i_g=\Ad u_g^i, u_g^iu_h^i=v^i_{g,h}u_{gh}^i$, for all $g$, $h\in G$. 
Let $\Cal{N}_{M_1\rtimes_{\alpha^1} G_1}(M_1)=
\{u\in\Cal{U}(M_1\rtimes_{\alpha^1} G_1)| uM_1u^*=M_1\}$ be the 
normalizer of $M_1$ in $M_1\rtimes_{\alpha^1} G_1$. 

There exists an isomorphism 
$\Cal{N}_{M_1\rtimes_{\alpha^1} G_1}(M_1)/\Cal{U}(M_1)\simeq G_1$ 
taking $u^1_g$ to $g$, and similarly 
$\Cal{N}_{M_2\rtimes_{\alpha^2} G_2}(M_2)/\Cal{U}(M_2)\simeq G_2$. 
Since $\Phi$ induces an isomorphism from 
$\Cal{N}_{M_1\rtimes_{\alpha^1} G_1}(M_1)$ to 
$\Cal{N}_{M_2\rtimes_{\alpha^2} G_2}(M_2)$, there exists a group 
isomorphism $\gamma:G_1 \to G_2$ such that
$$\Phi(u^1_g)=u^2_{\gamma(g)} (\text{mod }\Cal{U}(M_2))$$ 
for all $g\in G_1$.

Thus $\Phi(u^1_g)=w_gu^2_{\gamma(g)}$, for some unitaries 
$w_g\in \Cal{U}(M_2)$, $g \in G_1$. 
So $\Phi \alpha^1_g \Phi^{-1}=\Ad \Phi(u^1_g)=\Ad (w_gu^2_{\gamma(g)})=
\Ad w_g\text{ }\alpha^2_{\gamma(g)}$, for all $g\in G_1$, which proves (i).

Let $g$, $h \in G_1$. Then
$\Phi(v^1_{g,h})=\Phi(u^1_g)\Phi(u^1_h)\Phi(u^1_{gh})^*$
$= w_gu^2_{\gamma(g)}w_hu^2_{\gamma(h)}(u^2_{\gamma(gh)})^*\break w_{gh}^*$
$= w_g(\Ad u^2_{\gamma(g)})(w_h) 
u^2_{\gamma(g)}u^2_{\gamma(h)}(u^2_{\gamma(gh)})^*w_{gh}^*$ 
$= w_g\alpha^2_{\gamma(g)}(w_h)v^2_{\gamma(g),\gamma(h)}w^*_{gh}$, 
which proves (ii).

The converse follows easily by noticing that 
$\pi:L^2(M_1\rtimes_{\alpha^1} G_1)\simeq L^2(M_2\rtimes_{\alpha^2} G_2)$, 
defined by
$\pi(\sum_{g\in G} x_gu^1_g)=\sum_{g\in G} \Phi(x_g)w_gu^2_{\gamma(g)}$, 
for all $x_g\in M_1,g\in G$, is a Hilbert space 
isomorphism intertwining the $M_1\rtimes_{\alpha^1} G_1$ resp.
$M_2\rtimes_{\alpha^2} G_2$-actions (and hence the $M_1$ and $M_2$-actions)
on $L^2(M \rtimes_{\alpha^1} G_1)$ resp. $L^2(M \rtimes_{\alpha^2} G_2)$. 
\qed
\enddemo

Recall the simple fact that if the subfactors
$N \subset M$ and $\tilde N \subset \tilde M$ are isomorphic,
then the basic constructions $M_1$ and $\tilde M_1$ are isomorphic
as well (see e.g. \cite{Jo4}).

\proclaim{Corollary 1.8} Let $(\alpha^1,v^1)$, $(\alpha^2,v^2)$ be 
cocycle actions of the finite groups $G_1$, $G_2$ on the II$_1$ factors 
$M_1$, $M_2$. For $i=1$, $2$ let 
$M_i^{G_i}=\{x\in M_i| \alpha^i_g(x)=x$, for all 
$g\in G_i\}$. If there exists an isomorphism of inclusions
$$\Phi:(M_1^{G_1}\subset M_1) \to (M_2^{G_2} \subset M_2)$$
then $\alpha^1$ and $\alpha^2$ are outer conjugate via $\Phi$.
\endproclaim

\demo{Proof} Since $M_i^{G_i}\subset M_i\subset M_i\rtimes_{\alpha^i} G_i$,
$i=1,2$, is the basic construction (\cite{Jo4}), the isomorphism 
$\Phi$ can be extended to 
$M_1\rtimes_{\alpha^1} G_1\simeq M_2\rtimes_{\alpha^2} G_2$. 
Proposition 1.7 implies then that 
$\alpha^1$ and $\alpha^2$ are outer conjugate via $\Phi$. \qed
\enddemo

The next lemma shows that the ``restriction'' of an action of a group
$G$ on $M$ to the reduced algebra $pMp$, $p$ a projection in $M$, 
gives rise to a cocycle action of $G$ on $pMp$.

\proclaim{Lemma 1.9} Let $G$ be a discrete group, $\alpha$ an action of 
$G$ on the II$_1$ factor $M$, and $p$ a non-zero projection in $M$. 
Then there exists a cocycle action $\beta$ of $G$ on $pMp$ such that:

$$(pMp\subset p(M\rtimes_{\alpha} G)p)\simeq (pMp\subset pMp\rtimes_\beta G)$$

\endproclaim

\demo{Proof} Let $\tau$ be the unique faithful normalized trace of $M$. 
Since $\tau(\alpha_g(p))=\tau(p)$, for all $g\in G$, there exist 
unitaries $w_g\in \Cal U(M)$ such that $\alpha_g(p)=w_g^*pw_g$, 
for all $g\in G$. Let $\beta(g)=
\Ad w_g \alpha_g$. Denote by $u_g$ the unitary that implements $\alpha_g$ 
on $M$, for all $g\in G$. Thus $\beta_g$ is implemented by $w_gu_g$. By 
lemma 1.4, $\beta$ is a cocycle action of $G$ on $M$ with cocycle 
$v_{g,h}=w_g\alpha_g(w_h)w^*_{gh}$. Note that this is of course a 
coboundary for the $G$-action on $M$, but it may not be a coboundary 
when restricted to $pMp$. 

We show that $\beta$ is a cocycle action of $G$ on $pMp$ with cocycle $vp$. 
$\beta_g$ is an automorphism of $pMp$, since 
$\beta_g(p)=w_g\alpha_g(p)w_g^*=p$, 
for all $g\in G$. Applying to $p$ the relation $\beta_g\beta_h=
\Ad v_{g,h}\beta_{gh}$ implies that $p=\Ad v_{g,h}(p)$. Thus $p$ commutes 
with $v_{g,h}$, for all $g,h\in G$, so $v_{g,h}p$ is a cocycle for 
the cocycle action $\beta$ restricted to 
$pMp$. We have $v_{g,h}pv_{gh,k}p=\beta_g(v_{h,k}p)v_{g,hk}p$, for all 
$g$, $h$, $k\in G$.

Since $\alpha$, $\beta$ are cocycle conjugate, the inclusions 
$(M\subset M\rtimes_\alpha G)$ and $(M\subset M\rtimes_\beta G)$ 
are isomorphic, through an isomorphism that can be assumed to be the 
identity on $M$ (lemma 1.4, proposition 1.7). Hence we can identify 
$M\rtimes_\alpha G$ with $M\rtimes_\beta G$ (as von Neumann algebras), 
the unitaries implementing $\beta$ being identified with 
$w_gu_g\in M\rtimes_\alpha G$. The inclusions 
$(pMp\subset p(M\rtimes_{\alpha} G)p)$ and 
$(pMp\subset p(M\rtimes_{\beta} G)p)$ are isomorphic, so we only have 
to show that $p(M\rtimes_{\beta} G)p=pMp\rtimes_\beta G$ 
(we identify the abstract crossed product $pMp\rtimes_\beta G$ with 
its realization inside $M\rtimes_{\beta} G$). 

The von Neumann algebra $p(M\rtimes_{\beta} G)p$ is generated by 
$pxu_gp$ ($x\in M, g\in G$). Since $pxu_gp=(pxw_g^*p)(w_gu_g)$, it 
follows that $(pxu_gp)_{x\in M,g\in G}$ generate $pMp\rtimes_\beta G$ 
as a von Neumann algebra, so $p(M\rtimes_{\beta} G)p=pMp\rtimes_\beta G$.
\qed
\enddemo

\remark{Remark 1.10} It is easy to see that the cocycle conjugacy class 
of the cocycle action $\beta$ on $pMp$ depends only on $t=\tau(p)$
(\cite{Po6}). We will denote the cocycle conjugacy class of 
$(\beta,v,pMp)$ by $(\alpha^t,v^t,M^t)$, and call it 
{\it the amplification of $\alpha$ by t}. For values of $t$ greater 
than 1, define $\alpha^t$ to be the $t/n$-amplification of the 
action $id\otimes\alpha$ of $G$ on $M_n(\Bbb C)\otimes M$, for some 
$n\geq t$. Note that $\alpha^t$ is a properly outer cocycle action when 
$\alpha$ is properly outer (\cite{Po6}).
\endremark

\proclaim{Definition 1.11 (\cite{Po6})} Let $G$ be a discrete group 
with (cocycle)
action $\alpha$ on the II$_1$ factor $M$. The fundamental group of the 
action $\alpha$ is  

$$\Cal F(\alpha)=\{t>0 \, | \, \alpha^t\text{ is outer conjugate to }
\alpha\}.$$

Similarly we define

$$\Cal F^c(\alpha)=\{t>0 \, | \, \alpha^t\text{ is cocycle conjugate to }
\alpha\}.$$
\endproclaim
$\Cal F(\alpha)$ is an outer conjugacy invariant of $\alpha$, and
$\Cal F^c(\alpha)$ is a cocycle conjugacy invariant of $\alpha$.
Note that $\Cal F^c(\alpha) \subset \Cal F(\alpha)$ (see 1.3).

Let $G$ be an infinite discrete group, $\tau_0$ the normalized trace on 
$M_2(\Bbb C)$, and let $(R,\tau)=
\overline{\otimes_{g\in G} (M_2(\Bbb C),\tau_0)}^w$ be a copy of 
the hyperfinite II$_1$ factor. The {\it (non-commutative) Bernoulli 
$G$-action} on $R$ is the action 
$\sigma:G\rightarrow \text{Aut}(R)$ defined as 
$\sigma_g(\otimes_{h\in G}x_h)=\otimes_{h\in G} x'_h$, where 
$x'_h= x_{g^{-1}h}$, and $\{x_h\}_{h\in G}$ is such that all but 
finitely many $x_h$ are equal to 1. It is easy to
see (and well-known) that 
$\sigma$ is a properly outer, ergodic action. 

The following rigidity theorem from \cite{Po6} will provide the main
examples to which we will apply our construction in the next 
section.

\proclaim{Theorem 1.12} Let $G$ be a w-rigid group and $\sigma$ the
Bernoulli $G$-action on $R$. Then $\Cal F(\sigma)=\{1\}$.
\endproclaim

More generally, any of the Connes-St{\o}rmer Bernoulli $G$-actions 
with countable spectrum considered in \cite{Po6} have countable 
fundamental group. 

We will recall next the notion of {\it relative fundamental group}.
Let $N \subset M$ be an inclusion of II$_1$ factors. The {\it 
relative fundamental group} $\Cal F (N \subset M)$ is defined as
$$\Cal F (N \subset M) = \{ t > 0 \, | \, (N \subset M)^t
\ {\text{\rm is isomorphic to}} \ N \subset M \}$$
 $(N \subset M)^t$
denotes as usual the $t$-amplification of $N \subset M$
(see \cite{Po1}, \cite{Po2}).
Observe that $\Cal F (N \subset M)$ is clearly a multiplicative subgroup
of $\Bbb R_+^*$. If $N \subset M$ is stable, i.e. splits a common copy 
of the hyperfinite II$_1$ factor $R$, then clearly
$\Cal F (N \subset M) = \Bbb R_+^*$. This happens for instance if
the subfactor $N \subset M$ is constructed from an initial commuting
square by iterating the basic construction. See \cite{Bi1}, \cite{Bi2} 
for more on this.

If $\Cal F(N \subset M)$ is at most countable, then at most countably many
of the inclusions $pNp \subset pMp$, where $p$ runs through the set of
inequivalent projections in $N$, are isomorphic (as inclusions). Observe 
that the subfactors $pNp \subset pMp$ have all the same standard invariant.

Note that if $G$ is a discrete group with cocylc action $(\alpha, v)$
on the II$_1$ factor $M$, then 
$\Cal F^c (\alpha)$ $=\Cal F(M \subset M \rtimes_{\alpha,v} G)$
(this follows simply from the definitions).

The following rigidity result (\cite{Po7}, \cite{NPS}) 
provides further examples of actions having a fundamental group 
which is at most countable.

\medskip
\proclaim{Theorem 1.13} Let $M$ be a separable II$_1$
factor. Assume there exists a diffuse von Neumann subalgebra $B$ such 
that $B \subset M$ is a rigid inclusion (in the sense of \cite{Po7}) 
and $B'\cap M \subset B$. Then the fundamental group of $M$ is at most 
countable.
\endproclaim

\medskip
\remark{Remark 1.14} Let $G$ be a countable discrete ICC (infinite
conjugacy classes) group with property (T) and let $\alpha$ be an 
outer and ergodic action of $G$ on the hyperfinite
II$_1$ factor $R$. For instance, let $\alpha$ be the Bernoulli
$G$-action on $R$ described above. Let $M = R \rtimes_\alpha G$ and
let $u_g$, $g \in G$, be the unitaries in $M$ implementing $\alpha$.
Set $B = \{u_g \, | \, g \in G \}''$, and note that $B \simeq L(G)$. 
Then
$B\subset M$ is rigid and $L(G)' \cap M = R^G = \Bbb C$. Thus, by the 
above theorem, we deduce
that the fundamental group of $R\rtimes_\alpha G$ is at most countable. 
More generally, the same is true if $G$ is w-rigid and the action 
$\alpha$ is mixing. In particular, the relative fundamental group 
$\Cal F(R\subset R\rtimes_\alpha G)=\Cal F^c(\alpha)$ is at most
countable. 
If in addition $H^2(G,\Bbb T)$ is countable (hence automatically finite), 
the next lemma implies 
that $\Cal F(\alpha)$ is at most countable.
\endremark

\proclaim{Lemma 1.15} Let $G$ be a discrete group with countable second 
cohomology group \break $H^2(G,\Bbb T)$. Let $M$ be a II$_1$ factor, 
$I\subset \Bbb R$ 
an uncountable set and $(\alpha^i,v^i)_{i\in I}$ non cocycle conjugate 
cocycle actions of $G$ on $M$. Then $(\alpha^i,v^i)_{i\in I}$ are non 
outer conjugate modulo a countable set, i.e. 
$I(i_0)=\{i\in I, (\alpha^i,v^i)\text{ outer conjugate to }
(\alpha^{i_0},v^{i_0})\}$ is at most countable for each $i_0\in I$. 

In particular, given uncountably many outer conjugate actions of $G$ 
on $M$, uncountably many of these are actually cocycle conjugate 
actions.
\endproclaim

\demo{Proof} Assume by contradiction that $I(i_0)$ is uncountable for 
some $i_0$. According to lemma 1.2, for every $i$ there exists 
$\mu^i$ scalar 2-cocycle such that the actions $(\alpha^i,v^i)$ 
and $(\alpha^{i_0},\mu_iv^{i_0})$ are cocycle conjugate. Since 
$H^2(G,\Bbb T)$ is countable and $I(i_0)$ is uncountable, there exist 
$j_1,j_2\in I(i_0)$ such that $\mu^{j_1}\bar{\mu^{j_2}}$ is a coboundary. 
But $(\alpha^{j_1},v^{j_1})$ is cocycle conjugate to 
$(\alpha^{i_0},\mu^{j_1}v^{i_0})$, which is cocycle conjugate to 
$(\alpha^{j_2},\mu^{j_1}\bar{\mu^{j_2}}v^{j_2})$. Thus $\alpha^{j_1}$ 
and $\alpha^{j_2}$ are cocycle conjugate, which is a contradiction.
\enddemo

\bigskip

\heading
2. The construction
\endheading
\medskip

We consider in this section the class of subfactors introduced in 
\cite{BH}. Let $M$ be a II$_1$ factor and let $G$ be a countable
discrete group with an outer action $\alpha$ on $M$. Suppose
$G = \langle H, K \rangle$ is generated by two finite groups
$H$ and $K$. The subfactor $M^H \subset M \rtimes_\alpha K$
has index $|H|\cdot |K|$ and is
irreducible if and only if $H \cap K = \{ e \}$, where
$e$ denotes the identity in $G$. Note that we could start with
a cocycle action of $G$ and $M$. By \cite{Su1}, \cite{Su2}
we can modify the induced cocycle actions of the {\it finite}
groups $H$ and $K$ to actual actions. 

It is shown in \cite{BH}, \cite{BP} that many analytical and
algebraic properties of the subfactor $M^H \subset M \rtimes K$ 
are reflected by properties of the group $G$.
For instance, the following result is shown in \cite{BP}.

\medskip
\proclaim{Theorem 2.1} Let $H$ and $K$ be two finite groups with outer
actions $\sigma$ resp. $\rho$ on the II$_1$ factor $M$. Then the standard 
invariant of $M^H \subset M \rtimes K$ has property (T) (\cite{Po5}) 
if and only if 
the group $G$ generated by $\sigma (H)$ and $\rho (K)$ in the outer 
autmorphism group of $M$ has Kazhdan's property (T).
\endproclaim

\medskip
We will see below that the next theorem can be used to construct
continuous families of non-isomorphic, irreducible, finite index subfactors
of the hyperfinite II$_1$ factor all having the same standard invariant.
The construction can be carried out in such a way that this standard 
invariant will have property (T).

\medskip
The main result of this article is the following theorem.

\proclaim{Theorem 2.2}
Let $H$ be a finite abelian group and let $K = \Bbb Z_q$ be a cyclic 
group, where $q$ is a prime number. Let $G = \langle H , K \rangle$ be an 
infinite ICC group generated by $H$ and $K$.  Let $\alpha$ be a properly 
outer and ergodic action of $G$ on $R$. 

Then $\Cal F(R^H \subset R \rtimes_\alpha K)\subseteq \Cal F(\alpha)$. 
Hence, if $\Cal F(\alpha)$ is countable (resp. trivial), one obtains 
uncountably many (resp. a 1-parameter family of) irreducible 
subfactors of the hyperfinite II$_1$ factor $R$, which are 
non-isomorphic, but have all the same standard invariant.
\endproclaim
\medskip

Before we prove this theorem, let us give several examples of
groups and actions satisfying the hypothesis. 

If $G$ has property (T) with $H^2(G,\Bbb T)$ at most countable, 
and $\alpha$ is any properly outer and ergodic action of $G$ on $R$, 
then we established that $\Cal F(\alpha)$ is at most countable in
remark 1.14 and lemma 1.15. For instance, the groups 
$G_n= SL(2n+1, \Bbb Z)$ have Kazhdan's property
$(T)$ by \cite{K} and are ICC (see also \cite{HV}). They are 
(2,3)-generated for $n \ge 14$ by \cite{TWG} (see also \cite{TW}), 
i.e. $G_n$ is a quotient of the free product of 
$H= \Bbb Z_2$ and $K = \Bbb Z_3$ (this free product is of 
course just $PSL(2, \Bbb Z)$).
It follows from results of Steinberg (\cite{S1}, \cite{S2}, see 
also \cite{M}) that the second cohomology group 
$H^2(SL(n, \Bbb Z), \Bbb T)$ is a finite group (in fact it is equal 
to $\Bbb Z_2$) for $n \ge 5$.
These groups provide therefore (countably many) examples of groups 
satisfying the hypothesis of our theorem. We would like to thank 
Marsden Conder for
pointing out reference \cite{TW} and Pierre de la Harpe for the references
\cite{S1-2}, \cite{M}.

Recall that, by theorem 1.12, if $G$ is any w-rigid group and 
$\sigma$ the Bernoulli $G$-action on $R$, then 
$\Cal F(\sigma)=\{1\}$. Thus if $G=SL(2n+1,\Bbb Z),n\geq 14$, 
and $\alpha=\sigma$, we obtain one-parameter families of non-isomorphic, 
irreducible, index 6 hyperfinite subfactors having the same standard 
invariant. This standard invariant has property (T) (theorem 2.1).

A much larger class of examples can be obtained as follows: Let $G_1$ 
be the free product of any finite abelian group and a cyclic group of 
prime order. If $G_1\not= \Bbb Z_2 *\Bbb Z_2$ then $G_1$ is a hyperbolic 
group. Let $G_2$ be any hyperbolic property (T) group. By results of 
Olshanskii (see for instance \cite{Ol}, \cite{AMO}), there exists an infinite 
hyperbolic group $G$ which is a common quotient of $G_1$ and $G_2$. 
In particular, $G$ has property (T) and is generated by a finite 
abelian subgroup and a subgroup of prime order. Note that $G$ is ICC, 
since any non-elementary hyperbolic group is ICC. Thus, $G$ together 
with the Bernoulli $G$-action (or more generally any of the 
Connes-St{\o}rmer Bernoulli $G$-actions with countable spectrum 
considered in \cite{Po6}) satisfies the hypothesis of our theorem. 
We would like to thank Mark Sapir for pointing out this class of examples.

Note that similar results as in theorem 2.2 can be obtained by 
using the rigidity results in \cite{Po1} rather than theorems 1.12
and 1.13 quoted above.

\medskip
We proceed with the proof of theorem 2.2. We start with some lemmas.

\proclaim{Lemma 2.3} Let $H$ be a finite abelian group, let 
$K = \Bbb Z_q$ be a cyclic group, where $q$ is a prime number. Suppose
that $G = \langle H , \Bbb Z_q \rangle \ne H \cdot \Bbb Z_q$. Let
$\alpha$ be an outer cocycle action of $G$ on the hyperfinite 
II$_1$ factor $R$ and let $Q$ be the von Neumann algebra generated by 
the normalizer of $R^H$ in $R \rtimes_\alpha K$ 
(notation: $Q = \Cal N_{R \rtimes K}(R^H)''$). Then $Q = R$.
\endproclaim

\demo{Proof} 
We have by definition that $Q = \{ u \in \Cal U (R \rtimes K) \, | \,
u R^H u^* = R^H \} ''$. Since $H$ is abelian, 
we conclude that $(R^H \subset R) \cong (R_0 \subset R_0 \rtimes H)$, for 
some $R_0 \cong R$. Hence $R = \Cal N_R (R^H)''$ and we obtain therefore 
the chain of inclusions $R^H \subset R \subset Q \subset R \rtimes K $.
Since $R \subset R \rtimes K$ has no intermediate subfactors by
[Bi3, Theorem 3.2] we must have either $Q=R$ or $Q = R \rtimes K$.
If $Q = R \rtimes K$, then $(R^H \subset R \rtimes K) \cong
(R \subset R \rtimes (\Cal N(R^H) / \Cal U(R^H)))$ ([Jo1], [Jo3]).
Hence $R^H \subset R \rtimes K$ would have depth 2, contradicting the
fact that $G \ne H \cdot \Bbb Z_q$ ([BH]). Thus indeed $Q=R$. \qed
\enddemo

\medskip
\proclaim{Corollary 2.4} Let $H$ be a finite abelian group and
$K = \Bbb Z_q$, $q$ a prime number. Suppose that $G = \langle H, K \rangle$
is an infinite group and let $\alpha^i$ be outer cocycle actions of $G$ 
on the hyperfinite II$_1$ factor $R_i$, $i=1$, $2$.
Suppose that there is a surjective $*$-isomorphism 
$\Phi : R_1 \rtimes_{\alpha^1} K \to R_2 \rtimes_{\alpha^2}K$
such that $\Phi (R_1^H) = R_2^H$. Then $\Phi (R_1) = R_2$. In particular
we have $(R_1^H \subset R_1) \overset{\Phi} \to \cong (R_2^H \subset R_2)$ 
and $(R_1 \subset R_1 \rtimes_{\alpha^1} K ) \overset{\Phi} \to \cong
(R_2 \subset R_2 \rtimes_{\alpha^2}K)$.
\endproclaim

\demo{Proof} 
We have seen in Lemma 2.3 that $\Cal N_{R_i^H}(R_i \rtimes_{\alpha^i} K)''
= R_i$, $i=1$, $2$. But every (surjective) $*$-isomorphism takes 
normalizers to normalizers.  \qed
\enddemo

\proclaim{Proposition 2.5} Let $H$ be a finite abelian group and
$K = \Bbb Z_q$, $q$ a prime number. Suppose that $G = \langle H, K \rangle$
is an infinite group and let $\alpha^i$ be outer cocycle actions of $G$ on the hyperfinite II$_1$ factor $R_i$, $i=1$, $2$.
Suppose that there is a surjective $*$-isomorphism 
$\Phi : R_1 \rtimes_{\alpha^1} K \to R_2 \rtimes_{\alpha^2}K$
such that $\Phi (R_1^H) = R_2^H$. Then the cocycle actions $\alpha^1,\alpha^2$ of $G$ are outer conjugate by $\Phi$.
\endproclaim

\demo{Proof} It follows from corollary 2.4 and proposition 1.7 that 
$\Phi\alpha^1(K)\Phi^{-1}=\alpha^2(K)$ in Out($R_2$). From corollary 2.4 
and corollary 1.8 we deduce that 
$\Phi\alpha^1(H)\Phi^{-1}=\alpha^2(H)$ in Out($R_2$). 
Since $K$ and $H$ generate $G$, this implies 
$\Phi\alpha^1(G)\Phi^{-1}=\alpha^2(G)$ in Out($R_2$). \qed 
\enddemo

We give now the proof of theorem 2.2. 

\demo{Proof}
Let $G = \langle H , K \rangle $ be a quotient of the free product $H * K$
as in the theorem. Since $G$ is infinite, we have $G\not=H\cdot K$. 
Since $K$ is of prime order and $K\not\subset H$, it follows that 
$H\cap K=\{e\}$.

Let $\alpha$ be an outer and ergodic action of $G$
on the hyperfinite II$_1$ factor $R$. Let 
$t\in \Cal F(R^H\subset R\rtimes_\alpha K)$, $0 < t < 1$ (which
is sufficient since $\Cal F(R^H\subset R\rtimes_\alpha K)$ is a
group). We will show that $t\in \Cal F(\alpha)$.
 
Let $p$ be a projection in $R^H$ such that $\tau(p)=t$ ($\tau$ denotes 
as usual the normalized trace of $R$). Thus, the inclusions 
$R^H\subset R\rtimes_\alpha K$ and 
$pR^Hp\subset p(R\rtimes_\alpha K)p$ are isomorphic. 

By lemma 1.9, there exists a cocycle action ($\beta,v)$ of G on 
$pRp$ and an isomorphism 
$\Phi:(pRp\subset p(R\rtimes_\alpha G)p)\simeq (pRp\subset 
pRp\rtimes_{\beta} G)$, which is the identity on $pRp$. Moreover, 
from the construction of $\Phi$ (see lemma 1.9) it follows that 
$\Phi$ takes $p(R\rtimes_\alpha K)p$ onto $pRp\rtimes_{\beta} K$. 
Since $p\in R^H$, the actions $\alpha,\beta$ coincide on $H$ so the 
fixed point algebra $R^H$ is the same for both actions. Hence 
$\Phi$ takes $pR^Hp\subset pRp \subset p(R\times_\alpha K)p$ onto 
$pR^Hp\subset pRp\subset pRp\rtimes_{\beta} K$. Since $\beta_g(p)=p$, 
for all $g\in G$, we have $pR^Hp=(pRp)^H$ as subalgebras of 
$pRp\rtimes_{\beta} G$. This yields

$$(R^H\subset R\rtimes_\alpha K)\simeq (pR^Hp\subset p(R\rtimes_\alpha K)p)
\simeq ((pRp)^H\subset pRp \rtimes_{\beta} K)$$

Thus there exists an isomorphism $(R^H\subset R\rtimes_\alpha K)
\simeq ((pRp)^H\subset pRp \rtimes_{\beta} K)$. Proposition 2.5 implies 
that the cocycle actions $\alpha,\beta$ of $G$ are outer conjugate. 
Since $\beta$ is cocycle conjugate to $\alpha^t$, this implies 
$t\in \Cal F(\alpha)$ which ends the proof.  \qed
\enddemo

\enddocument
 
\vfill
\bye